\documentclass[10pt,twoside]{article}
\usepackage{amsfonts}
\usepackage{fancyhdr}
\usepackage{titlesec}
\usepackage{cite}
\usepackage{ifthen}
\usepackage{pifont}
\usepackage{stmaryrd}
\usepackage{setspace}
\usepackage{indentfirst}
\usepackage{amsmath,amssymb,amscd,bbm,amsthm,mathrsfs,dsfont}
\input amssym.def

\titleformat{\section}{\centering\large\bfseries}{\S\arabic{section}}{1em}{}
\newboolean{first}
\setboolean{first}{true}

\newtheorem{theorem}{Theorem}[section]
\newtheorem{lemma}{Lemma}[section]
\newtheorem{rem}{Remark}[section]

\newtheorem{definition}{Definition}[section]
\begin{document}

\setlength\abovedisplayskip{2pt}
\setlength\abovedisplayshortskip{0pt}
\setlength\belowdisplayskip{2pt}
\setlength\belowdisplayshortskip{0pt}

\title{\bf \Large The Boundedness of some Multilinear operators with rough kernel on the weighted Morrey spaces\author{HE Sha\quad  TAO Xiang-xing}\date{}} \maketitle
 \footnote{Received: 2012-10-8.}
 \footnote{MR Subject Classification: 42B20, 42B25, 42B35.}
 \footnote{Digital Object Identifier(DOI): 10.1007/s11766-012-****-*.}
 \footnote{Partially supported by the National Natural Science Foundation of China\,(11171306, 11071065)
 and sponsored by the Scientific Project of Zhejiang Provincial Science Technology Department\,(2011C33012) 
 and the Scientific Research Fund of Zhejiang Provincial Education Department\,(Z201017584).}
\begin{center}
\begin{minipage}{135mm}

{\bf \small Abstract}.\hskip 2mm {\small
In this paper, strong boundedness of $T_{\Omega,\alpha}^A$ and $M_{\Omega,\alpha}^A$, the multilinear fractional integral operators and the corresponding fractional maximal operators, are showed on weighted Morrey spaces with two weights when $D^{\gamma}A \in \dot{\Lambda_{\beta}}(|\gamma|=m-1)$ or $D^{\gamma}A \in BMO(|\gamma|=m-1)$. For the mulitilinear singular integral operators $T_{\Omega}^A$ and the corresponding maximal operators $M_{\Omega}^A$, they are proved to be strong bounded operators on the same spaces if $D^{\gamma}A \in \dot{\Lambda_{\beta}}(|\gamma|=m-1)$; and if $D^{\gamma}A \in BMO(|\gamma|=m-1)(m=1, 2)$, the boundedness of $T_{\Omega}^A$ and $M_{\Omega}^A$ are obtained on weighted Morrey spaces with one weight.}
\end{minipage}
\end{center}

\thispagestyle{fancyplain} \fancyhead{}
\fancyhead[L]{\textit{Appl. Math. J. Chinese Univ.}\\
2012, 27(*): ***-***} \fancyfoot{} \vskip 10mm

\section{Introduction}
Let us consider the following multilinear fractional integral operator:
\begin{equation*}
T_{\Omega,\alpha}^{A}f(x)=\int_{\mathbb{R}^n}\frac{\Omega(x-y)}{|x-y|^{n-\alpha+m-1}}R_m(A;x,y)f(y)dy\qquad0<\alpha<n
\end{equation*}
and the corresponding multilinear fractional maximal operator:
\begin{equation*}
M_{\Omega,\alpha}^{A}f(x)=\sup\limits_{r>0}\frac{1}{r^{n-\alpha+m-1}}\int_{|x-y|<r}|\Omega(x-y)R_m(A;x,y)f(y)|dy\qquad0<\alpha<n
\end{equation*}
where $\Omega\in L^s(S^{n-1})(s>1)$ is homogeneous of degree zero in $\mathbb{R}^n$, $A$ is a function defined on $\mathbb{R}^n$ and $R_m(A;x,y)$ denotes the $m$-th order Taylor series remainder of $A$ at $x$ expanded about $y$, that is,
\begin{equation*}
R_m(A;x,y)=A(x)-\sum\limits_{|\gamma|<m}\frac{1}{\gamma!}D^{\gamma}A(y)(x-y)^{\gamma}
\end{equation*}
$\gamma=(\gamma_1,\cdots,\gamma_n)$, each $\gamma_i(i=1,\cdots,n)$ is a nonnegative integer, $|\gamma|=\sum\limits_{i=1}^{n}\gamma_i$, $\gamma!=\gamma_1!\cdots\gamma_n!$, $x^{\gamma}=x_1^{\gamma_1}\cdots x_n^{\gamma_n}$ and $D^{\gamma}=\frac{\partial^{|\gamma|}}{\partial^{\gamma_1}x_1\cdots\partial^{\gamma_n}x_n}$.

We notice that when $\alpha=0$, the above operators become the multilinear singular integral operator and the corresponding maximal operator:
\begin{equation*}
T_{\Omega}^{A}f(x)=\int_{\mathbb{R}^n}\frac{\Omega(x-y)}{|x-y|^{n+m-1}}R_m(A;x,y)f(y)dy
\end{equation*}
\begin{equation*}
M_{\Omega}^{A}f(x)=\sup\limits_{r>0}\frac{1}{r^{n+m-1}}\int_{|x-y|<r}|\Omega(x-y)R_m(A;x,y)f(y)|dy
\end{equation*}
For $m=1$, $T_{\Omega,\alpha}^{A}$ is obviously the commutator operator, $[A,T_{\Omega,\alpha}]f(x)=A(x)T_{\Omega,\alpha}f(x)-T_{\Omega,\alpha}(Af)(x)$, where $T_{\Omega,\alpha}$ is the following fractional integral operator:
\begin{equation*}
T_{\Omega,\alpha}f(x)=\int_{\mathbb{R}^n}\frac{\Omega(x-y)}{|x-y|^{n-\alpha}}f(y)dy\qquad0<\alpha<n
\end{equation*}

The classical Morrey spaces were first introduced by Morrey to study the local behavior of solutions to second order elliptic partial differential equations. From then on, a lot of works concerning Morrey spaces and some related spaces have been done, see [21], [23-28] for details. In 2009, Komori and Shirai [14] considered the weighted Morrey spaces and investigated some classical singular integrals in harmonic analysis on them, such as the Hardy-Littlewood maximal operator, the Calder\'{o}n-Zygmund operator, the fractional integral operator as well as the fractional maximal operator and so on. Recently, Wang [30] discussed the weighted boundedness of the classical singular operators with rough kernels on the weighted Morrey spaces.

In recent years, the multilinear theorey have attracted much attentions. In 1967, Bajsanski and Coifman [1] proved the boundedness of the multilinear operator associated with the commutators of singular integrals considered by Calder\'{o}n. In [2], Cohen studied the $L^p$ boundedness of the operator $T_{\Omega}^A$ for $m=2$. Using the method of good-$\lambda$ inequality, in 1986, Cohen and Gosselin [3] proved that if $\Omega\in Lip_1(S^{n-1})$,  satisfies vanishing condition and $D^{\gamma}A\in BMO(|\gamma|=m-1)$, then $T_{\Omega}^A$ is bounded on $L^p$. In 1994, for $m=2$, Hofmann [12] proved that $T_{\Omega}^A$ is a bounded operator on $L^p(w)$ when $\Omega\in L^{\infty}(S^{n-1})$ and $w\in A_p$. In 2001, Ding obtained $T_{\Omega,\alpha}^A$ and $M_{\Omega,\alpha}^A$ are both weighted bounded operators from $L^p(w^p)$ to $L^q(w^q)$ with $w\in A(p,q)$ and from $L^p(1\leq p<\frac{n}{\alpha})$ to $L^{\frac{n}{n-\alpha},\infty}$ with power weight when $D^{\gamma}A\in L^r(1<r\leq \infty, |\gamma|=m-1)$ in [5]. Later, Ding and Lu [8] proved the $(L^p(w^p), L^q(w^q))$ boundedness of $T_{\Omega,\alpha}^{A_1,\cdots,A_k}$ and its corresponding maximal operator $M_{\Omega,\alpha}^{A_1,\cdots,A_k}$ (the definition of them will be given later). In 2002, Wu and Yang [31] studied if $D^{\gamma}A\in BMO(|\gamma|=m-1)$, then $T_{\Omega,\alpha}^A$ is bounded on $L^p$. After that, Lu and Zhang [17] obtained $T_{\Omega,\alpha}^A$ is a bounded operator from $L^p$ to $L^q$ $(1<p<\frac{n}{\alpha+\beta})$ and from $L^1$ to $L^{\frac{n}{n-(\alpha+\beta)},\infty}$ when $D^{\gamma}A\in \dot{\Lambda}_{\beta}(|\gamma|=m-1)$. In the same year, Lu, Wu and Zhang [16] proved $T_{\Omega}^A$ and $M_{\Omega}^A$ is bounded from $L^p$ to $L^q$ $(\frac{1}{q}=\frac{1}{p}-\frac{\beta}{n})$ and from $L^1$ to $L^{\frac{n}{n-\beta},\infty}$ with the bound $C\sum\limits_{|\gamma|=m-1}\|D^{\gamma}A\|_{\dot{\Lambda}_{\beta}}$. In [13], for $m=2$ Jiao showed $T_{\Omega}^A$ is bounded on $L^p$ and bounded from $H^1$ to weak $L^1$ and from $L^{\infty}$ to $BMO$ if $\Omega$ satisfies some conditions. In 2005, Ding [6] obtained the two-weight boundedness of $T_{\Omega,\alpha}^A$ under the condition $D^{\gamma}A\in L^r(|\gamma|=m-1)$. In 2008, Han and Lu [11] obtained the $(H^1,L^{{\frac{n}{n-\alpha}},\infty})$ boundedness of $T_{\Omega,\alpha}^A$ when $D^{\gamma}A\in BMO(|\gamma|=m-1)$. Recently, Wu and Tao [32] studied if $D^{\gamma}A_i\in BMO(|\gamma_i|=m_i-1)$, then $T_{\Omega,\alpha}^{A_1,\cdots,A_k}$ is bounded from $H^1$ to $L^{{\frac{n}{n-\alpha}},\infty}$.

From the above review, we find that many works concerning $T_{\Omega,\alpha}^A$, $M_{\Omega,\alpha}^A$, $T_{\Omega}^A$ and $M_{\Omega}^A$ have been done on $L^p$ spaces when $D^{\gamma}A$ belongs to $L^p$, $BMO$ or Lipschitz spaces if $|\gamma|=m-1$, but there are not any study about these operators on weighted Morrey spaces. So our purpose in this paper is to consider the above operators on weighted Morrey spaces.

The organization of this paper is as follows. We will introduce in next section some definitions and notations that are necessary. The main results will be given in Section 3. In Section 4, we give some requisite lemmas and well-known results that are important in proving theorems. The proof of the theorems will be shown in Section 5.
\section{Definitions and notations}
A weight is a locally integrable function on $\mathbb{R}^n$ which takes values in $(0,\infty)$ almost everywhere. For a weight $w$ and a measurable set $E$, we define $w(E)=\int_{E}w(x)dx$, the Lebesgue measure of $E$ by $|E|$ and the characteristic function of $E$ by $\chi_{E}$. The weighted Lebesgue spaces with respect to the measure $w(x)dx$ are denoted by $L^p(w)$ with $0<p<\infty$. we say a weight $w$ satisfies the doubling condition if there exists a constant $D>0$ such that for any ball $B$, we have $w(2B)\leq Dw(B)$. When $w$ satisfies this condition, we denote $w\in \Delta_2$ for short.

Troughout this paper, $B(x_0,r)$ denotes a ball centered at $x_0$ with radius $r$. $Q$ be a cube with sides parallel to the axes. For $K>0$, $KQ$ denotes the cube with the same center as $Q$ and side-length being $K$ times longer. When $\alpha=0$, we will denote $T_{\Omega,\alpha}$, $T_{\Omega,\alpha}^A$, $M_{\Omega,\alpha}^A$ by $T_{\Omega}$, $T_{\Omega}^A$, $M_{\Omega}^A$ respectively. And for any number $a$, $a'$ is standing for the conjugate of $a$. The letter $C$ is used for various constants, and may changes from one occurrence to another.

To begin with, we introduce the weighted Morrey spaces.
\begin{definition}{[14]}
Let $1\leq p<\infty$, $0<\kappa<1$ and $w$ be a weight. Then a weighted Morrey space is defined by
\begin{equation*}
L^{p,\kappa}(w):=\{f\in L_{loc}^{p}(w):\|f\|_{L^{p,\kappa}(w)}<\infty\}
\end{equation*}
where
\begin{equation*}
\|f\|_{L^{p,\kappa}(w)}=\sup\limits_{B}\Big(\frac{1}{w(B)^{\kappa}}\int_{B}|f(x)|^pw(x)dx\Big)^{\frac{1}{p}}
\end{equation*}
and the supremum is taken over all balls $B$ in $\mathbb{R}^n$.
\end{definition}
In the case of fractional order, we need to consider a weighted Morrey space with two weights. It is defined as follows:
\begin{definition}{[14]}
Let $1\leq p<\infty$, $0<\kappa<1$, $u$, $v$ be two weights. The two weights weighted Morrey space is defined by
\begin{equation*}
L^{p,\kappa}(u,v):=\{f:\|f\|_{L^{p,\kappa}(u,v)}<\infty\}
\end{equation*}
where
\begin{equation*}
\|f\|_{L^{p,\kappa}(u,v)}=\sup\limits_{B}\Big(\frac{1}{v(B)^{\kappa}}\int_{B}|f(x)|^pu(x)dx\Big)^{\frac{1}{p}}
\end{equation*}
and the supremum is taken over all balls $B$ in $\mathbb{R}^n$. If $u=v$, then we denote $L^{p,\kappa}(u)$ for short.
\end{definition}
From Remark 2.2 in [14], we could define the weighted Morrey spaces with cubes instead of balls. So we shall use these two definitions of weighted Morrey spaces appropriate to calculation.

Next, we give the definition of Lipschitz space and $BMO$ space.
\begin{definition}\label{def2.3}
The Lipschitz space of order $\beta$, $0<\beta<1$ is defined by
\begin{equation*}
\dot{\Lambda}_\beta(\mathbb{R}^n)=\{f:|f(x)-f(y)|\leq C|x-y|^{\beta}\}
\end{equation*}
and the smallest constant $C>0$ is the Lipschitz norm $\|\cdot\| _{\dot{\Lambda}_{\beta}}$.
\end{definition}
\begin{definition}
A locally integrable function $b$ is said to be in $BMO(\mathbb{R}^n)$ if
\begin{equation*}
\|b\|_*=\|b\|_{BMO}=\sup\limits_B\frac{1}{|B|}\int_B|b(x)-b_B|dx<\infty
\end{equation*}
where
\begin{equation*}
b_B=\frac{1}{|B|}\int_Bb(y)dy
\end{equation*}
\end{definition}

At last, we shall show the definition of two weight classes.
\begin{definition}
A weight function $w$ is in the Muckenhoupt class $A_p$ with $1<p<\infty$ if there exists $C>1$ such that for any ball B
\begin{equation*}
\Big(\frac{1}{|B|}\int_{B}w(x)dx\Big)\Big(\frac{1}{|B|}\int_{B}w(x)^{-\frac{1}{p-1}}dx\Big)^{p-1}\leq C
\end{equation*}
We define $A_{\infty}=\bigcup\limits_{1<p<\infty}A_p$.\\
When $p=1$, $w\in A_1$ if there exists $C>1$ such that for almost every $x$,
\begin{equation*}
Mw(x)\leq Cw(x)
\end{equation*}
\end{definition}
\begin{definition}{[18]}
A weight function $w$ belongs to $A(p,q)$ for $1<p<q<\infty$ if there exists $C>1$ such that
\begin{equation*}
\Big(\frac{1}{|B|}\int_Bw(x)^qdx\Big)^{\frac{1}{q}}\Big(\frac{1}{|B|}\int_Bw(x)^{-\frac{p}{p-1}}dx\Big)^{\frac{p-1}{p}}\leq C
\end{equation*}
When $p=1$, $w$ is in $A(1,q)$ with $1<q<\infty$ if there exists $C>1$ such that
\begin{equation*}
\Big(\frac{1}{|B|}\int_Bw(x)^qdx\Big)^{\frac{1}{q}}\Big(ess\sup\limits_{x\in B}\frac{1}{w(x)}\Big)\leq C
\end{equation*}
\end{definition}
\begin{rem}{[14]}\label{rem2.7}\quad
If $w\in A(p,q)$ with $1<p<q$, then \\
(a)$w^q$, $w^{-p'}$, $w^{-q'}\in \Delta_2$.\\
(b)$w^{-p'}\in A_{t'}$ with $1/t+1/t'=1$.
\end{rem}

\section{Main Theorems}
Our main results will be stated as follows.
\begin{theorem}\label{th3.1}
If $0<\alpha+\beta<n$, $\Omega\in L^s(S^{n-1})(s>1)$ is homogeneous of degree zero, $1<s'<p<n/(\alpha+\beta)$, $1/q=1/p-(\alpha+\beta)/n$, $0<\kappa<p/q$, $w^{s'}\in A(p/s',q/s')$, $D^{\gamma}A\in \dot{\Lambda}_{\beta}(|\gamma|=m-1)$, then
\begin{align}\label{eq3.1}
\|T_{\Omega,\alpha}^{A}f\|_{L^{q,\kappa q/p}(w^q)}\leq C\sum\limits_{|\gamma|=m-1}\|D^{\gamma}A\|_{\dot{\Lambda}_{\beta}}\|f\|_{L^{p,\kappa}(w^p,w^q)}
\end{align}
\begin{align}\label{eq3.2}
\|M_{\Omega,\alpha}^{A}f\|_{L^{q,\kappa q/p}(w^q)}\leq C\sum\limits_{|\gamma|=m-1}\|D^{\gamma}A\|_{\dot{\Lambda}_{\beta}}\|f\|_{L^{p,\kappa}(w^p,w^q)}
\end{align}
\end{theorem}
\begin{theorem}\label{th3.2}
If $0<\beta<1$, $\Omega\in L^s(S^{n-1})(s>1)$ is homogeneous of degree zero, $1<s'<p<n/\beta$, $1/q=1/p-\beta/n$, $0<\kappa<p/q$, $w^{s'}\in A(p/s',q/s')$, $D^{\gamma}A\in \dot{\Lambda}_{\beta}(|\gamma|=m-1)$, then
\begin{align}\label{eq3.3}
\|T_{\Omega}^{A}f\|_{L^{q,\kappa q/p}(w^q)}\leq C\sum\limits_{|\gamma|=m-1}\|D^{\gamma}A\|_{\dot{\Lambda}_{\beta}}\|f\|_{L^{p,\kappa}(w^p,w^q)}
\end{align}
\begin{align}\label{eq3.4}
\|M_{\Omega}^{A}f\|_{L^{q,\kappa q/p}(w^q)}\leq C\sum\limits_{|\gamma|=m-1}\|D^{\gamma}A\|_{\dot{\Lambda}_{\beta}}\|f\|_{L^{p,\kappa}(w^p,w^q)}
\end{align}
\end{theorem}
\begin{theorem}\label{th3.3}
If $0<\alpha<n$, $\Omega\in L^s(S^{n-1})(s>1)$ is homogeneous of degree zero, $1<s'<p<n/\alpha$, $1/q=1/p-\alpha/n$, $0<\kappa<p/q$, $w^{s'}\in A(p/s',q/s')$, $D^{\gamma}A\in BMO(|\gamma|=m-1)$, then
\begin{align}\label{eq3.5}
\|T_{\Omega,\alpha}^Af\|_{L^{q,\kappa q/p}(w^q)}\leq C\sum\limits_{|\gamma|=m-1}\|D^{\gamma}A\|_*\|f\|_{L^{p,\kappa}(w^p,w^q)}
\end{align}
\begin{align}\label{eq3.6}
\|M_{\Omega,\alpha}^Af\|_{L^{q,\kappa q/p}(w^q)}\leq C\sum\limits_{|\gamma|=m-1}\|D^{\gamma}A\|_*\|f\|_{L^{p,\kappa}(w^p,w^q)}
\end{align}
\end{theorem}
For $T_{\Omega}^A$ and $M_{\Omega}^A$, we only study the cases when $m=1$ and $m=2$. In these two cases, we denote $T_{\Omega}^A$, $M_{\Omega}^A$ by $[A,T_{\Omega}]$, $[A,M_{\Omega}]$ and $\tilde{T}_{\Omega}^A$, $\tilde{M}_{\Omega}^A$ respectively in order to distinguish from $T_{\Omega}^A$ and $M_{\Omega}^A$ that are suitable for any integer $m$. That is,
\begin{align*}
&[A,T_{\Omega}]f(x)=\int_{\mathbb{R}^n}\frac{\Omega(x-y)}{|x-y|^{n}}\big(A(x)-A(y)\big)f(y)dy\\
&[A,M_{\Omega}]f(x)=\sup\limits_{r>0}\frac{1}{r^{n}}\int_{|x-y|<r}\Omega(x-y)\big(A(x)-A(y)\big)f(y)dy
\end{align*}
and
\begin{align*}
&\tilde{T}_{\Omega}^Af(x)=\int_{\mathbb{R}^n}\frac{\Omega(x-y)}{|x-y|^{n+1}}\big(A(x)-A(y)-\nabla A(y)(x-y)\big)f(y)dy\\
&\tilde{M}_{\Omega}^Af(x)=\sup\limits_{r>0}\frac{1}{r^{n+1}}\int_{|x-y|<r}\Omega(x-y)\big(A(x)-A(y)-\nabla A(y)(x-y)\big)f(y)dy
\end{align*}
Then for the above operators, we have the following results on weighted Morrey spaces with one weight.
\begin{theorem}\label{th3.4}
If $\Omega\in L^{s}(S^{n-1})(s>1)$ is homogeneous of degree zero and satisfies the vanishing condition $\int_{S^{n-1}}\Omega(x')d\sigma(x')=0$, $1<s'<p<\infty$, $0<\kappa<1$, $w\in A_{p/s'}$, $A\in BMO$, then
\begin{align}\label{eq3.7}
\|[A,T_{\Omega}]f\|_{L^{p,\kappa}(w)}\leq C\|A\|_*\|f\|_{L^{p,\kappa}(w)}
\end{align}
\begin{align}\label{eq3.8}
\|[A,M_{\Omega}]f\|_{L^{p,\kappa}(w)}\leq C\|A\|_*\|f\|_{L^{p,\kappa}(w)}
\end{align}
\end{theorem}
\begin{theorem}\label{th3.5}
If $\Omega\in L^{\infty}(S^{n-1})$ is homogeneous of degree zero and satisfies the moment condition $\int_{S^{n-1}}\theta\Omega(\theta)d\theta=0$, $1<p<\infty$, $0<\kappa<1$, $w\in A_p$, $\nabla A\in BMO$, then
\begin{align}\label{eq3.9}
\|\tilde{T}_{\Omega}^Af\|_{L^{p,\kappa}(w)}\leq C\|\nabla A\|_*\|f\|_{L^{p,\kappa}(w)}
\end{align}
\begin{align}\label{eq3.10}
\|\tilde{M}_{\Omega}^Af\|_{L^{p,\kappa}(w)}\leq C\|\nabla A\|_*\|f\|_{L^{p,\kappa}(w)}
\end{align}
\end{theorem}
Here we point out that for $T_{\Omega}^A$ and $M_{\Omega}^A$, when $D^{\gamma}A\in BMO(|\gamma|=m-1)(m\geq 3)$, the analogous conclusions of Theorem \ref{th3.5} is open.

{\bf{Remark}}\qquad
Define
\begin{align*}
&T_{\Omega,\alpha}^{A_1,\cdots,A_k}f(x)=\int_{\mathbb{R}^n}\prod\limits_{i=1}^kR_{m_i}(A_i;x,y)\frac{\Omega(x-y)}{|x-y|^{n-\alpha+N}}f(y)dy\\
&M_{\Omega,\alpha}^{A_1,\cdots,A_k}f(x)=\sup\limits_{r>0}\frac{1}{r^{n-\alpha+N}}\int_{|x-y|<r}|\Omega(x-y)|\prod\limits_{i=1}^k|R_{m_i}(A_i;x,y)||f(y)|dy
\end{align*}
when $0<\alpha<n$, it becomes a class of multilinear fractional integral operators. When $\alpha=0$, it is a class of multilinear singular integral operators. $R_{m_i}(A_i;x,y)=A_i(x)-\sum\limits_{|\gamma|<m_i}\frac{1}{\gamma!}D^{\gamma}A_i(y)(x-y)^{\gamma}\quad(i=1,\cdots,k)$, $N=\sum\limits_{i=1}^k(m_i-1)$. Repeating the proofs of theorems above, we will find that, for $T_{\Omega,\alpha}^{A_1,\cdots,A_k}$ and $M_{\Omega,\alpha}^{A_1,\cdots,A_k}$ the conclusions of Theorem 3.1 and Theorem 3.2 above with the bounds $C\prod\limits_{i=1}^k\big(\sum\limits_{|\gamma|=m_i-1}\|D^{\gamma}A_i\|_{\dot{\Lambda}_{\beta}}\big)$ and Theorem 3.3 with the bounds $C\prod\limits_{i=1}^k\big(\sum\limits_{|\gamma|=m_i-1}\|D^{\gamma}A_i\|_*\big)$ also hold, respectively.

\section{Lemmas and well-known results}
\begin{lemma}{[3]}\label{lem1}
Let $A(x)$ be a function on $\mathbb{R}^n$ with $m$-th order derivatives in $L_{loc}^l(\mathbb{R}^n)$ for some $l>n$. Then
\begin{align*}
|R_m(A;x,y)|\leq C|x-y|^m\sum\limits_{|\gamma|=m}\Big(\frac{1}{|I_x^y|}\int_{I_x^y}|D^{\gamma}A(z)|^ldz\Big)^{\frac{1}{l}}
\end{align*}
where $I_x^y$ is the cube centered at $x$ with sides parallel to the axes,whose diameter is $5\sqrt{n}|x-y|$.
\end{lemma}
\begin{lemma}{[20]}\label{lem2}
For $0<\beta<1$, $1\leq q<\infty$, we have
\begin{align*}
\|f\|_{\dot{\Lambda}_{\beta}}\approx \sup\limits_Q\frac{1}{|Q|^{1+\beta/n}}\int_Q|f(x)-m_Q(f)|dx\approx \sup\limits_Q\frac{1}{|Q|^{\beta/n}}\Big(\frac{1}{|Q|}\int_Q|f(x)-m_Q(f)|^qdx\Big)^{\frac{1}{q}}
\end{align*}
Here $m_Q(f)$ denotes the average of $f$ over $Q$.\\
For $q=\infty$, the formula should be interpreted appropriately.
\end{lemma}
\begin{lemma}{[9]}\label{lem3}
Let $Q_1\subset Q_2$, $g\in \dot{\Lambda}_{\beta}(0<\beta<1)$. Then
\begin{align*}
|m_{Q_1}(g)-m_{Q_2}(g)|\leq C|Q_2|^{\beta/n}\|g\|_{\dot{\Lambda}_{\beta}}
\end{align*}
\end{lemma}
\begin{theorem}{[7]}\label{lem4}
Suppose that $0<\alpha<n$, $1<p<n/\alpha$, $1/q=1/p-\alpha/n$ and $\Omega\in L^s(S^{n-1})(s>1)$. Then $T_{\Omega,\alpha}$ is a bound operator from $L^p(w^p)$ to $L^q(w^q)$, if $s,p,q$ and $w$ satisfy one of the following conditions.\\
(a)$1\leq s'<p$ and $w(x)^{s'}\in A(p/s',q/s')$\\
(b)$s>q$ and $w(x)^{-s'}\in A(q'/s',p'/s')$\\
(c)$s>1$, $\alpha/n+1/s<1/p<1/s'$ and there is an $r$ s.t. $1<r<s/(n/\alpha)'$ and $w(x)^{r'}\in A(p,q)$.
\end{theorem}
\begin{lemma}{[14]}\label{lem5}
If $w\in \Delta_2$, then there exists a constant $D_1>1$, s.t.
\begin{align}\label{eq4.5}
w(2B)\geq D_1w(B).
\end{align}
\end{lemma}
We call $D_1$ the reverse doubling constant.
\begin{theorem}{[8]}\label{lem6}
Suppose that $0<\alpha<n$, $1<p<n/\alpha$, $1/q=1/p-\alpha/n$, $\Omega$ is homogeneous of degree zero with $\Omega\in L^s(S^{n-1})(s>1)$. Moreover, for $1\leq j\leq k$, $|\gamma|=m_j-1$, $m_j\geq 2$, and $D^{\gamma}A_j\in BMO(\mathbb{R}^n)$. If the index set $\{\alpha, p, q, s\}$ satisfies one of the following conditions:\\
(a) $s'<p$, $w(x)^{s'}\in A(p/s',q/s')$;\\
(b) $s>q$, $w(x)^{-s'}\in A(q'/s',p'/s')$;\\
(c) $\alpha/n+1/s<1/p<1/s'$, there is an $r$, $1<r<s/(n/\alpha)'$ such that $w(x)^{r'}\in A(p,q)$;\\
then there is a $C>0$, independent of $f$ and $A_j$, such that
\begin{align*}
\Big(\int_{\mathbb{R}^n}|T_{\Omega,\alpha}^{A_1,\cdots,A_k}f(x)w(x)|^qdx\Big)^{1/q}\leq C\prod\limits_{j=1}^k\Big(\sum\limits_{|\gamma|=m_j-1}\|D^{\gamma}A_j\|_*\Big)\Big(\int_{\mathbb{R}^n}|f(x)w(x)|^pdx\Big)^{1/p}
\end{align*}
\end{theorem}
\begin{lemma}{[4]}\label{lem7}
(John-Nirenberg Lemma) Let $1\leq p<\infty$. Then $b\in BMO$ if and only if
\begin{align*}
\frac{1}{|Q|}\int_Q|b-b_Q|^pdx\leq C\|b\|_*^p
\end{align*}
\end{lemma}
\begin{lemma}{[22]}\label{lem8}
Assume $b\in BMO$, then for cube $Q_1\subset Q_2$,
\begin{align*}
|b_{Q_1}-b_{Q_2}|\leq Clog(|Q_2|/|Q_1|)\|b\|_*
\end{align*}
\end{lemma}
From the above Lemma, we immediately get the following conclusion:
\begin{lemma}{[29]}\label{lem9}
If $b\in BMO$, then
\begin{align*}
|b_{2^{j+1}B}-b_{B}|\leq 2^n(j+1)\|b\|_*
\end{align*}
\end{lemma}
\begin{theorem}{[15]}\label{lem10}
Suppose that $\Omega\in L^s(S^{n-1})(s>1)$ is homogeneous of degree zero and satisfies the vanishing condition $\int_{S^{n-1}}\Omega(x')d\sigma(x')=0$. If $b\in BMO(\mathbb{R}^n)$ and $p,q,w$ satisfy one of the following conditions, then $[b,T_\Omega]$ is bounded on $L^p(w)$:
(a) $s'\leq p<\infty$, $p\neq 1$ and $w\in A_{p/s'}$;\\
(b) $1 \leq p\leq s$, $p\neq \infty$ and $w^{1-p'}\in A_{p'/s'}$;\\
(c) $1\leq p<\infty$ and $w^{s'}\in A_p$.
\end{theorem}
\begin{theorem}{[12]}\label{lem11}
If $\Omega\in L^{\infty}(S^{n-1})$ is homogeneous of degree zero and satisfies the moment condition $\int_{S^{n-1}}\theta\Omega(\theta)d\theta=0$, $w\in A_p$, $1<p<\infty$, $\nabla A\in BMO$, then we have
\begin{align*}
\|\tilde{T}_{\Omega}^Af\|_{L^p(w)}\leq C\|\Omega\|_{\infty}\|\nabla A\|_*\|f\|_{L^p(w)}
\end{align*}
\end{theorem}
\begin{lemma}{[10]}\label{lem12}
The following are true:\\
(1) If $w\in A_p$ for some $1\leq p<\infty$, then $w\in \Delta_2$. More precisely, for all $\lambda >1$ we have
\begin{align*}
w(\lambda Q)\leq C\lambda ^{np}w(Q).
\end{align*}
(2) If $w\in A_p$ for some $1\leq p<\infty$, then there exist $C>0$ and $\delta>0$ such that for any cube $Q$ and a measurable set $S\subset Q$,
\begin{align*}
\frac{w(S)}{w(Q)}\leq C\Big(\frac{|S|}{|Q|}\Big)^\delta.
\end{align*}
\end{lemma}
\begin{lemma}{[19]}\label{lem13}
Let $w\in A_{\infty}$. Then the norm of $BMO(w)$ is equivalent to the norm of $BMO(\mathbb{R}^n)$, where
\begin{align*}
&BMO(w)=\Big\{b:\|b\|_{*,w}=\sup\limits_Q\frac{1}{w(Q)}\int_{Q}|b(x)-m_{Q,w}b|w(x)dx\Big\}\\
and\\
&m_{Q,w}b=\frac{1}{w(Q)}\int_Qb(x)w(x)dx
\end{align*}
\end{lemma}

\section{Proofs of the Main Results}
{\bf{Proof of Theorem \ref{th3.1}}}.\\
To prove (\ref{eq3.1}), we give a pointwise estimat of $T_{\Omega,\alpha}^Af(x)$ at first. Set
\begin{align*}
\bar{T}_{\Omega,\alpha+\beta}f(x)=\int_{\mathbb{R}^n}\frac{|\Omega(x-y)|}{|x-y|^{n-\alpha-\beta}}|f(y)|dy\qquad0<\alpha+\beta<n
\end{align*}
where $\Omega\in L^s(S^{n-1})(s>1)$ is homogeneous of degree zero in $\mathbb{R}^n$.
Then we have the following theorem:
\begin{theorem}\label{th5.1}
If $\alpha\geq 0$, $0<\alpha+\beta<n$, $D^{\gamma}A\in \dot{\Lambda}_{\beta}(|\gamma|=m-1)$, then there exists a constant $C$ independent of $f$ such that
\begin{align*}
|T_{\Omega,\alpha}^Af(x)|\leq C\Big(\sum\limits_{|\gamma|=m-1}\|D^{\gamma}A\|_{\dot{\Lambda}_{\beta}}\Big)\bar{T}_{\Omega,\alpha+\beta}f(x)
\end{align*}
\end{theorem}
{\bf{Proof}}\quad For fixed $x\in \mathbb{R}^n$, $r>0$, let $Q$ be a cube centered at $x$ and has diameter $r$, $Q_k=2^kQ$ and set
\begin{align}\label{eq5.2}
A_{Q_k}(y)=A(y)-\sum\limits_{|\gamma|=m-1}\frac{1}{\gamma!}m_{Q_k}(D^{\gamma}A)y^{\gamma}
\end{align}
where $m_{Q_k}f$ is the average of $f$ on $Q_k$. Then we have when $|\gamma|=m-1$
\begin{align}\label{eq5.3}
D^{\gamma}A_{Q_k}(y)=D^{\gamma}A(y)-m_{Q_k}(D^{\gamma}A)
\end{align}
and from [3], we have
\begin{align}\label{eq5.4}
R_m(A;x,y)=R_m(A_{Q_k};x,y)
\end{align}
Then,
\begin{align*}
|T_{\Omega,\alpha}^Af(x)|\leq \sum\limits_{k=-\infty}^{\infty}\int_{2^{k-1}r\leq |x-y|<2^kr}\frac{|R_m(A_{Q_k};x,y)|}{|x-y|^{m-1}}\frac{|\Omega(x-y)|}{|x-y|^{n-\alpha}}|f(y)|dy:=\sum\limits_{k=-\infty}^{\infty}T_k
\end{align*}
By Lemma \ref{lem1} we get
\begin{align*}
&|R_m(A_{Q_k};x,y)|\leq |R_{m-1}(A_{Q_k};x,y)|+C\sum\limits_{|\gamma|=m-1}|D^{\gamma}A_{Q_k}(y)||x-y|^{m-1}\\
&\leq C|x-y|^{m-1}\sum\limits_{|\gamma|=m-1}\Big(\frac{1}{|I_x^y|}\int_{I_x^y}|D^{\gamma}A_{Q_k}(z)|^ldz\Big)^{\frac{1}{l}}+C|x-y|^{m-1}\sum\limits_{|\gamma|=m-1}|D^{\gamma}A_{Q_k}(y)|
\end{align*}
Note that, if $|x-y|<2^kr$, then $I_x^y\subset 5nQ_k$. By Lemma \ref{lem2} and Lemma \ref{lem3} we have when $|\gamma|=m-1$,
\begin{align*}
&\Big(\frac{1}{|I_x^y|}\int_{I_x^y}|D^{\gamma}A_{Q_k}(z)|^ldz\Big)^{\frac{1}{l}}=\Big(\frac{1}{|I_x^y|}\int_{I_x^y}|D^{\gamma}A(z)-m_{Q_k}(D^{\gamma}A)|^ldz\Big)^{\frac{1}{l}}\\
&\leq \Big(\frac{1}{|I_x^y|}\int_{I_x^y}|D^{\gamma}A(z)-m_{I_x^y}(D^{\gamma}A)|^ldz\Big)^{\frac{1}{l}}+|m_{I_x^y}(D^{\gamma}A)-m_{5nQ_k}(D^{\gamma}A)|\\
&\qquad+|m_{5nQ_k}(D^{\gamma}A)-m_{Q_k}(D^{\gamma}A)|\\
&\leq C|Q_k|^{\frac{\beta}{n}}\|D^{\gamma}A\|_{\dot{\Lambda}_{\beta}}\leq C(2^kr)^{\beta}\|D^{\gamma}A\|_{\dot{\Lambda}_{\beta}}
\end{align*}
From Definition \ref{def2.3}, we obtain when $|\gamma|=m-1$,
\begin{align*}
&|D^{\gamma}A_{Q_k}(y)|=|D^{\gamma}A(y)-m_{Q_k}(D^{\gamma}A)|
\leq C|Q_k|^{\frac{\beta}{n}}\|D^{\gamma}A\|_{\dot{\Lambda}_{\beta}}
\leq C(2^kr)^{\beta}\|D^{\gamma}A\|_{\dot{\Lambda}_{\beta}}
\end{align*}
Thus, \begin{align*}
|R_m(A_{Q_k};x,y)|\leq C|x-y|^{m-1}(2^kr)^{\beta}\sum\limits_{|\gamma|=m-1}\|D^{\gamma}A\|_{\dot{\Lambda}_{\beta}}
\end{align*}
Therefore,
\begin{align*}
&T_k\leq C\sum\limits_{|\gamma|=m-1}\|D^{\gamma}A\|_{\dot{\Lambda}_{\beta}}\int_{2^{k-1}r\leq |x-y|<2^kr}\frac{(2^kr)^{\beta}}{|x-y|^{n-\alpha}}|\Omega(x-y)||f(y)|dy\\
&\leq C\sum\limits_{|\gamma|=m-1}\|D^{\gamma}A\|_{\dot{\Lambda}_{\beta}}\int_{2^{k-1}r\leq |x-y|<2^kr}\frac{|\Omega(x-y)|}{|x-y|^{n-\alpha-\beta}}|f(y)|dy
\end{align*}
It follows that
\begin{align*}
&|T_{\Omega,\alpha}^Af(x)|\leq \sum\limits_{k=-\infty}^{\infty}\Big(C\sum\limits_{|\gamma|=m-1}\|D^{\gamma}A\|_{\dot{\Lambda}_{\beta}}\int_{2^{k-1}r\leq |x-y|<2^kr}\frac{|\Omega(x-y)|}{|x-y|^{n-\alpha-\beta}}|f(y)|dy\Big)\\
&\leq C\sum\limits_{|\gamma|=m-1}\|D^{\gamma}A\|_{\dot{\Lambda}_{\beta}}\sum\limits_{k=-\infty}^{\infty}\int_{2^{k-1}r\leq |x-y|<2^kr}\frac{|\Omega(x-y)|}{|x-y|^{n-\alpha-\beta}}|f(y)|dy\\
&=C\sum\limits_{|\gamma|=m-1}\|D^{\gamma}A\|_{\dot{\Lambda}_{\beta}}\int_{\mathbb{R}^n}\frac{|\Omega(x-y)|}{|x-y|^{n-\alpha-\beta}}|f(y)|dy\\
&=C\sum\limits_{|\gamma|=m-1}\|D^{\gamma}A\|_{\dot{\Lambda}_{\beta}}\bar{T}_{\Omega,\alpha+\beta}f(x)
\end{align*}
Thus we finish the proof of theorem \ref{th5.1}.

The following Theorem is a key theorem in proving (\ref{eq3.1}).
\begin{theorem}\label{th5.2}
Under the same conditions of Theorem \ref{th3.1}, $\bar{T}_{\Omega,\alpha+\beta}$ is bounded from $L^{p,\kappa}(w^p,w^q)$ to $L^{q,\kappa q/p}(w^q)$.
\end{theorem}
{\bf{Proof}}\quad Fix a ball $B(x_0,r_B)$ and We decompose $f=f_1+f_2$ with $f_1=f\chi_{2B}$. Since $\bar{T}_{\Omega,\alpha}$ is a linear operator, we have
\begin{align*}
&\|\bar{T}_{\Omega,\alpha+\beta}f\|_{L^{q,\kappa q/p}(w^q)}=
\Big(\frac{1}{w^q(B)^{\kappa q/p}}\int_B|\bar{T}_{\Omega,\alpha+\beta}f(x)|^qw^q(x)dx\Big)^{\frac{1}{q}}\\
&\leq \frac{1}{w^q(B)^{\kappa/p}}\Big(\int_B|\bar{T}_{\Omega,\alpha+\beta}f_1(x)|^qw^q(x)dx\Big)^{\frac{1}{q}}+
\frac{1}{w^q(B)^{\kappa/p}}\Big(\int_B|\bar{T}_{\Omega,\alpha+\beta}f_2(x)|^qw^q(x)dx\Big)^{\frac{1}{q}}\\
&=J_1+J_2
\end{align*}
We estimate $J_1$ at first. From Remark \ref{rem2.7} (a) we know that $w^q\in \Delta_2$, then by Theorem \ref{lem4} (a) we get,
\begin{align*}
&J_1\leq \frac{1}{w^q(B)^{\kappa/p}}\|\bar{T}_{\Omega,\alpha+\beta}f_1\|_{L^q(w^q)}\\
&\leq \frac{C}{w^q(B)^{\kappa/p}}\|f_1\|_{L^p(w^p)}
=\frac{C}{w^q(B)^{\kappa/p}}\Big(\int_{2B}|f(x)|^pw(x)^pdx\Big)^{\frac{1}{p}}\\
&\leq C\|f\|_{L^{p,\kappa}(w^p,w^q)}\frac{w^q(2B)^{\kappa/p}}{w^q(B)^{\kappa/p}}\leq C\|f\|_{L^{p,\kappa}(w^p,w^q)}
\end{align*}
Now we consider the term $J_2$.
\begin{align*}
&|\bar{T}_{\Omega,\alpha+\beta}f_2(x)|=\int_{(2B)^c}\frac{|\Omega(x-y)|}{|x-y|^{n-\alpha-\beta}}|f(y)|dy=\sum\limits_{j=1}^\infty\int_{2^{j+1}B\setminus2^jB}\frac{|\Omega(x-y)|}{|x-y|^{n-\alpha-\beta}}|f(y)|dy\\
&\leq C\sum\limits_{j=1}^\infty\big(\int_{2^{j+1}B}|\Omega(x-y)|^sdy\big)^{\frac{1}{s}}
\Big(\int_{2^{j+1}B\setminus2^jB}\frac{|f(y)|^{s'}}{|x-y|^{(n-\alpha-\beta)s'}}dy\Big)^{\frac{1}{s'}}=C\sum\limits_{j=1}^\infty(I_{1j}I_{2j})
\end{align*}
We will estimate $I_{1j}$, $I_{2j}$ respectively.
Let $z=x-y$, then for $x\in B$, $y\in 2^{j+1}B$, we have $z\in 2^{j+2}B$. Noticing that $\Omega$ is homogeneous of degree zero and $\Omega\in L^{s}(S^{n-1})$, we obtain
\begin{align*}
&I_{1j}=\big(\int_{2^{j+2}B}|\Omega(z)|^sdz\big)^{\frac{1}{s}}=\Big(\int_{0}^{2^{j+2}r_B}\int_{S^{n-1}}|\Omega(z')|^sdz'r^{n-1}dr\Big)^{\frac{1}{s}}\\
&=C\|\Omega\|_{L^{s}(S^{n-1})}|2^{j+2}B|^{\frac{1}{s}}
\end{align*}
where $z'=z/|z|$.
For $x\in B$, $y\in (2B)^c$, $|x-y|\sim|x_0-y|$, thus
\begin{align*}
&I_{2j}\leq\frac{C}{|2^{j+1}B|^{1-\frac{\alpha+\beta}{n}}}\big(\int_{2^{j+1}B}|f(y)|^{s'}dy\big)^{\frac{1}{s'}}
\end{align*}
By H$\ddot{o}$lder inequality and $w^{s'}\in A(p/s',q/s')$, we get
\begin{align*}
&\big(\int_{2^{j+1}B}|f(y)|^{s'}dy\big)^{\frac{1}{s'}}
\leq C\big(\int_{2^{j+1}B}|f(y)|^pw(y)^pdy\big)^{\frac{1}{p}}\big(\int_{2^{j+1}B}w(y)^{-\frac{ps'}{p-s'}}dy\big)^{{\frac{p-s'}{ps'}}}\\
&\leq C\|f\|_{L^{p,\kappa}(w^p,w^q)}w^q(2^{j+1}B)^{\frac{\kappa}{p}}\big(\int_{2^{j+1}B}w(y)^{-\frac{ps'}{p-s'}}dy\big)^{{\frac{p-s'}{ps'}}}\\
&\leq C\|f\|_{L^{p,\kappa}(w^p,w^q)}w^q(2^{j+1}B)^{\frac{\kappa}{p}}\frac{|2^{j+1}B|^{\frac{pq-s'q+s'p}{pqs'}}}{w^q(2^{j+1}B)^{\frac{1}{q}}}
\end{align*}
Thus
\begin{align*}
&|\bar{T}_{\Omega,\alpha+\beta}f_2(x)|\leq C\sum\limits_{j=1}^\infty(I_{1j}I_{2j})\leq C\sum\limits_{j=1}^\infty\|f\|_{L^{p,\kappa}(w^p,w^q)}\frac{1}{w^q(2^{j+1}B)^{\frac{1}{q}-\frac{\kappa}{p}}}
\end{align*}
So we get
\begin{align*}
J_2\leq C\|f\|_{L^{p,\kappa}(w^p,w^q)}\sum\limits_{j=1}^\infty\frac{w^q(B)^{\frac{1}{q}-\frac{\kappa}{p}}}{w^q(2^{j+1}B)^{\frac{1}{q}-\frac{\kappa}{p}}}
\end{align*}
Using Remark \ref{rem2.7} (a) and Lemma \ref{lem5}, we get that $w$ satisfies inequality (\ref{eq4.5}), so the above series converges since the reverse doubling constant is larger than one, as a result,
\begin{align*}
J_2\leq C\|f\|_{L^{p,\kappa}(w^p,w^q)}
\end{align*}
Therefore, we have showed the proof of Theorem \ref{th5.2}.\\
Here we remark that Theorem \ref{th5.2} is essentially verifying the multilinear fractional operator $T_{\Omega,\alpha}$ is bounded on weighted Morrey spaces.

Now let us turn to {\bf{prove inequality (\ref{eq3.1})}}. By Theorem \ref{th5.1} and Theorem \ref{th5.2}, (\ref{eq3.1}) immediately obtained.\\
We are now in a place of {\bf{proving (\ref{eq3.2}) in Theorem \ref{th3.1}}}.\\ Set
\begin{align*}
\bar{T}_{\Omega,\alpha}^Af(x)=\int_{\mathbb{R}^n}\frac{|\Omega(x-y)|}{|x-y|^{n-\alpha+m-1}}|R_m(A;x,y)||f(y)|dy\qquad0\leq\alpha<n
\end{align*}
where $\Omega\in L^s(S^{n-1})(s>1)$ is homogeneous of degree zero in $\mathbb{R}^n$.
It is easy to see that, for $\bar{T}_{\Omega,\alpha}^A$, the conclusions of inequality (\ref{eq3.1}) also hold. On the other hand, for any $r>0$, we have
\begin{align*}
&\bar{T}_{\Omega,\alpha}^Af(x)\geq\int_{|x-y|<r}\frac{|\Omega(x-y)|}{|x-y|^{n-\alpha+m-1}}|R_m(A;x,y)||f(y)|dy\\
&\geq\frac{1}{r^{n-\alpha+m-1}}\int_{|x-y|<r}|\Omega(x-y)||R_m(A;x,y)||f(y)|dy
\end{align*}
Taking the supremum for $r>0$ on the inequality above, we get
\begin{align}\label{eq5.1}
\bar{T}_{\Omega,\alpha}^Af(x)\geq M_{\Omega,\alpha}^Af(x)
\end{align}
Thus, we can immediately obtain (\ref{eq3.2}) from (\ref{eq5.1}) and (\ref{eq3.1}).\\

Before starting proving Theorem \ref{th3.2}, we give the following theorem at first since this theorem plays an important role in proving Theorem \ref{th3.2}. Set
\begin{align*}
\bar{T}_{\Omega,\beta}f(x)=\int_{\mathbb{R}^n}\frac{|\Omega(x-y)|}{|x-y|^{n-\beta}}f(y)dy
\end{align*}
where $\Omega\in L^s(S^{n-1})(s>1)$ is homogeneous of degree zero in $\mathbb{R}^n$.
\begin{theorem}\label{th5.3}
Under the assumptions of Theorem \ref{th3.2}, $\bar{T}_{\Omega,\beta}$ is bounded from $L^{p,\kappa}(w^p,w^q)$ to $L^{q,\kappa q/p}(w^q)$.
\end{theorem}
We shall omit the proof for it is the same as that of Theorem \ref{th5.2}.\\
Now, let us {\bf{prove Theorem \ref{th3.2}}}. It is not difficult to see that inequality (\ref{eq3.3}) of Theorem \ref{th3.2} can be easily obtained from Theorem \ref{th5.1} and Theorem \ref{th5.3}. At the same time, We can immediately arrive at (\ref{eq3.4}) from (\ref{eq5.1}) and (\ref{eq3.3}).\\
From now on, we are in the place of {\bf{showing Theorem \ref{th3.3}.}} We study (\ref{eq3.5}) at first.
Fixing any cube $Q$ whose center is $x$ and diameter is $r$, denote $\bar{Q}=2Q$ and set
\begin{align*}
A_{\bar{Q}}(y)=A(y)-\sum\limits_{|\gamma|=m-1}\frac{1}{\gamma!}m_{\bar{Q}}(D^{\gamma}A)y^{\gamma}
\end{align*}
we notice that the above equality is the special case of equality (\ref{eq5.2}) when $k=1$. Thus the equality (\ref{eq5.3}) and (\ref{eq5.4}) also hold for $A_{\bar{Q}}(y)$.
We decompose $f$ according to $\bar{Q}$, that is $f=f\chi_{\bar{Q}}+f\chi_{(\bar{Q})^c}:=f_1+f_2$. Then we have
\begin{align*}
&\|T_{\Omega,\alpha}^Af\|_{L^{q,\kappa q/p}(w^q)}
\leq \frac{1}{w^q(Q)^{\kappa/p}}\Big(\int_Q|T_{\Omega,\alpha}^Af_1(y)|^qw(y)^qdy\Big)^{\frac{1}{q}}\\&+
\frac{1}{w^q(Q)^{\kappa/p}}\Big(\int_Q|T_{\Omega,\alpha}^Af_2(y)|^qw(y)^qdy\Big)^{\frac{1}{q}}
=I+II
\end{align*}
From Remark \ref{rem2.7} (a), we know that $w^q\in \Delta_2$, then togethering with Theorem \ref{lem6} (a) we have
\begin{align*}
&I\leq \frac{C}{w^q(Q)^{\kappa/p}}\sum\limits_{|\gamma|=m-1}\|D^{\gamma}A\|_*(\int_{\bar{Q}}|f(y)|^pw(y)^pdy)^{\frac{1}{p}}\\
&=C \sum\limits_{|\gamma|=m-1}\|D^{\gamma}A\|_*\|f\|_{L^{p,\kappa}(w^p,w^q)}\Big(\frac{w^q(\bar{Q})}{w^q(Q)}\Big)^{\kappa/p}\\
&\leq C\sum\limits_{|\gamma|=m-1}\|D^{\gamma}A\|_*\|f\|_{L^{p,\kappa}(w^p,w^q)}
\end{align*}
Next, we consider the term $T_{\Omega,\alpha}^Af_2(y)$ contained in $II$. By Lemma \ref{lem1} and equality (\ref{eq5.3}), (\ref{eq5.4}), we have
\begin{align*}
&|T_{\Omega,\alpha}^Af_2(y)|\leq \int_{(\bar{Q})^c}\frac{|R_m(A_{\bar{Q}};y,z)|}{|y-z|^{m-1}}\frac{|\Omega(y-z)|}{|y-z|^{n-\alpha}}|f(z)|dz\\
&\leq C\int_{(\bar{Q})^c}\sum\limits_{|\gamma|=m-1}\Big(\frac{1}{|I_y^z|}\int_{I_y^z}|D^{\gamma}A_{\bar{Q}}(t)|^ldt\Big)^{\frac{1}{l}}|\Omega(y-z)|\frac{|f(z)|}{|y-z|^{n-\alpha}}dz\\
&\qquad+C\int_{(\bar{Q})^c}\sum\limits_{|\gamma|=m-1}|D^{\gamma}A(z)-m_{\bar{Q}}(D^{\gamma}A)||\Omega(y-z)|\frac{|f(z)|}{|y-z|^{n-\alpha}}dz=II_1+II_2
\end{align*}
We estimate $II_1$ and $II_2$ respectively.
By Lemma \ref{lem7} and Lemma \ref{lem8} and then use the similar steps as the proof of Theorem \ref{th5.2}, we get
\begin{align*}
II_1\leq C\sum\limits_{|\gamma|=m-1}\|D^{\gamma}A\|_*\|f\|_{L^{p,\kappa}(w^p,w^q)}\sum\limits_{j=1}^{\infty}\frac{1}{w^q(2^{j+1}Q)^{\frac{1}{q}-\frac{\kappa}{p}}}
\end{align*}
For $y\in Q$, $z\in (\bar{Q})^c$, we have $|y-z|\sim|x-z|$, so we obtain
\begin{align*}
&II_2\leq C\sum\limits_{j=1}^{\infty}\int_{2^{j+1}Q\setminus2^jQ}\sum\limits_{|\gamma|=m-1}|D^{\gamma}A(z)-m_{2^{j+1}Q}(D^{\gamma}A)||\Omega(y-z)|\frac{|f(z)|}{|y-z|^{n-\alpha}}dz\\
&+C\sum\limits_{j=1}^{\infty}\int_{2^{j+1}Q\setminus2^jQ}\sum\limits_{|\gamma|=m-1}|m_{2^{j+1}Q}(D^{\gamma}A)-m_{\bar{Q}}(D^{\gamma}A)||\Omega(y-z)|\frac{|f(z)|}{|y-z|^{n-\alpha}}dz\\
&=II_{21}+II_{22}
\end{align*}
By H$\ddot{o}$lder inequality, we get
\begin{align*}
&II_{21}\leq C\sum\limits_{j=1}^{\infty}\frac{1}{|2^jQ|^{1-\alpha/n}}\int_{2^{j+1}Q}\sum\limits_{|\gamma|=m-1}|D^{\gamma}A(z)-m_{2^{j+1}Q}(D^{\gamma}A)||\Omega(y-z)||f(z)|dz\\
&\leq C\sum\limits_{j=1}^{\infty}\frac{1}{|2^jQ|^{1-\alpha/n}}\Big(\int_{2^{j+1}Q}\sum\limits_{|\gamma|=m-1}|D^{\gamma}A(z)-m_{2^{j+1}Q}(D^{\gamma}A)|^{s'}|f(z)|^{s'}dz\Big)^{\frac{1}{s'}}\cdot\\
&\qquad\Big(\int_{2^{j+1}Q}|\Omega(y-z)|^sdz\Big)^{\frac{1}{s}}\\
&\leq C\|\Omega\|_{L^s(S^{n-1})}\sum\limits_{j=1}^{\infty}\frac{|2^{j+2}Q|^{\frac{1}{s}}}{|2^jQ|^{1-\alpha/n}}\Big(\int_{2^{j+1}Q}\sum\limits_{|\gamma|=m-1}|D^{\gamma}A(z)-m_{2^{j+1}Q}(D^{\gamma}A)|^{s'}|f(z)|^{s'}dz\Big)^{\frac{1}{s'}}\\
\end{align*}
\begin{align*}
&\leq C\|\Omega\|_{L^s(S^{n-1})}\sum\limits_{j=1}^{\infty}\frac{|2^{j+2}Q|^{\frac{1}{s}}}{|2^jQ|^{1-\alpha/n}}\big(\int_{2^{j+1}Q}|f(z)|^pw(z)^pdz\big)^{\frac{1}{p}}\cdot\\
&\qquad\Big(\int_{2^{j+1}Q}\sum\limits_{|\gamma|=m-1}|D^{\gamma}A(z)-m_{2^{j+1}Q}(D^{\gamma}A)|^{\frac{ps'}{p-s'}}w(z)^{-\frac{ps'}{p-s'}}dz\Big)^{\frac{p-s'}{ps'}}\\
\end{align*}
We can calculate the part of including the function $D^{\gamma}A$ as follows:
\begin{align*}
&\Big(\int_{2^{j+1}Q}\sum\limits_{|\gamma|=m-1}|D^{\gamma}A(z)-m_{2^{j+1}Q}(D^{\gamma}A)|^{\frac{ps'}{p-s'}}w(z)^{-\frac{ps'}{p-s'}}dz\Big)^{\frac{p-s'}{ps'}}\\
&\leq C\Big(\int_{2^{j+1}Q}\sum\limits_{|\gamma|=m-1}\Big|D^{\gamma}A(z)-m_{2^{j+1}Q,w^{-\frac{ps'}{p-s'}}}(D^{\gamma}A)\Big|^{\frac{ps'}{p-s'}}w(z)^{-\frac{ps'}{p-s'}}dz\Big)^{\frac{p-s'}{ps'}}\\
&\qquad+\sum\limits_{|\gamma|=m-1}\Big|m_{2^{j+1}Q,w^{-\frac{ps'}{p-s'}}}(D^{\gamma}A)-m_{2^{j+1}Q}(D^{\gamma}A)\Big|w^{-\frac{ps'}{p-s'}}(2^{j+1}Q)^{\frac{p-s'}{ps'}}=III+IV
\end{align*}
For the term $III$, as $w^{s'}\in A(p/s',q/s')$, using Remark \ref{rem2.7} (b), we obtain that $w^{-\frac{ps'}{p-s'}}\in A_{t'}\subset A_{\infty}$ ($1/t+1/t'=1$), then by Lemma \ref{lem13} that the norm of $BMO(w^{-\frac{ps'}{p-s'}})$ is equivalent to the norm of $BMO(\mathbb{R}^n)$ and $w^{s'}\in A(p/s',q/s')$ condition, we obtain
\begin{align*}
&III\leq C\sum\limits_{|\gamma|=m-1}\|D^{\gamma}A\|_*w^{-\frac{ps'}{p-s'}}(2^{j+1}Q)^{\frac{p-s'}{ps'}}\\
&=C\sum\limits_{|\gamma|=m-1}\|D^{\gamma}A\|_*\frac{|2^{j+1}Q|^{{\frac{p-s'}{ps'}}+\frac{1}{q}}}{w^q(2^{j+1}Q)^{\frac{1}{q}}}
\end{align*}
For the term $IV$, from the John and Nirenberg lemma, we have that there exist $C_1>0$ and $C_2>0$ such that for any cube $Q$ and $s>0$
\begin{align*}
\Big|\Big\{t\in 2^{j+1}Q:\sum\limits_{|\gamma|=m-1}|D^{\gamma}A(t)-m_{2^{j+1}Q}(D^{\gamma}A)|>s\Big\}\Big|\leq C_1|2^{j+1}Q|e^{{-C_2s/\big(\sum\limits_{|\gamma|=m-1}\|D^{\gamma}A\|_*}\big)}
\end{align*}
since $\sum\limits_{|\gamma|=m-1}(D^{\gamma}A)\in BMO$. Then from Lemma \ref{lem12} (2), we have
\begin{align*}
w\Big(\Big\{t\in 2^{j+1}Q:\sum\limits_{|\gamma|=m-1}|D^{\gamma}A(t)-m_{2^{j+1}Q}(D^{\gamma}A)|>s\Big\}\Big)\leq Cw(2^{j+1}Q)e^{{-C_2s\delta/\big(\sum\limits_{|\gamma|=m-1}\|D^{\gamma}A\|_*}\big)}
\end{align*}
for some $\delta>0$. Hence the inequality implies
\begin{align*}
&\sum\limits_{|\gamma|=m-1}\Big|m_{2^{j+1}Q,w^{-\frac{ps'}{p-s'}}}(D^{\gamma}A)-m_{2^{j+1}Q}(D^{\gamma}A)\Big|\\
&\leq \frac{1}{w^{-\frac{ps'}{p-s'}}(2^{j+1}Q)}\int_{2^{j+1}Q}\sum\limits_{|\gamma|=m-1}|D^{\gamma}A(t)-m_{2^{j+1}Q}(D^{\gamma}A)|w^{-\frac{ps'}{p-s'}}(t)dt\\
&=C\frac{1}{w^{-\frac{ps'}{p-s'}}(2^{j+1}Q)}\int_0^{\infty}w^{-\frac{ps'}{p-s'}}(\{t\in 2^{j+1}Q:\sum\limits_{|\gamma|=m-1}|D^{\gamma}A(t)-m_{2^{j+1}Q}(D^{\gamma}A)|>s\})ds
\end{align*}
\begin{align*}
&\leq C\frac{1}{w^{-\frac{ps'}{p-s'}}(2^{j+1}Q)}\int_0^{\infty}w^{-\frac{ps'}{p-s'}}(2^{j+1}Q)e^{{-C_2s\delta/\big(\sum\limits_{|\gamma|=m-1}\|D^{\gamma}A\|_*}\big)}ds
=C\sum\limits_{|\gamma|=m-1}\|D^{\gamma}A\|_*
\end{align*}
So,
\begin{align*}
&IV\leq C\sum\limits_{|\gamma|=m-1}\|D^{\gamma}A\|_*w^{-\frac{ps'}{p-s'}}(2^{j+1}Q)^{\frac{p-s'}{ps'}}\\
&=C\sum\limits_{|\gamma|=m-1}\|D^{\gamma}A\|_*\frac{|2^{j+1}Q|^{\frac{p-s'}{ps'}+\frac{1}{q}}}{w^q(2^{j+1}Q)^{\frac{1}{q}}}
\end{align*}
Thus,
\begin{align*}
II_{21}\leq C\sum\limits_{|\gamma|=m-1}\|D^{\gamma}A\|_*\|f\|_{L^{p,\kappa}(w^p,w^q)}\sum\limits_{j=1}^{\infty}\frac{1}{w^q(2^{j+1}Q)^{\frac{1}{q}-\frac{\kappa}{p}}}
\end{align*}
For the term $II_{22}$, using Lemma \ref{lem9} and the analogous steps as the proof of Theorem \ref{th5.2}, we get
\begin{align*}
&II_{22}\leq C\sum\limits_{|\gamma|=m-1}\|D^{\gamma}A\|_*\sum\limits_{j=1}^{\infty}j\int_{2^{j+1}Q\setminus2^jQ}|\Omega(y-z)|\frac{|f(z)|}{|y-z|^{n-\alpha}}dz\\
&\leq C\sum\limits_{|\gamma|=m-1}\|D^{\gamma}A\|_*\|f\|_{L^{p,\kappa}(w^p,w^q)}\sum\limits_{j=1}^{\infty}\frac{j}{w^q(2^{j+1}Q)^{\frac{1}{q}-\frac{\kappa}{p}}}
\end{align*}
Therefore,
\begin{align*}
|T_{\Omega,\alpha}^Af_2(y)|\leq C\sum\limits_{|\gamma|=m-1}\|D^{\gamma}A\|_*\|f\|_{L^{p,\kappa}(w^p,w^q)}\sum\limits_{j=1}^{\infty}\frac{j}{w^q(2^{j+1}Q)^{\frac{1}{q}-\frac{\kappa}{p}}}
\end{align*}
Consequently,
\begin{align*}
&II\leq C\sum\limits_{|\gamma|=m-1}\|D^{\gamma}A\|_*\|f\|_{L^{p,\kappa}(w^p,w^q)}\sum\limits_{j=1}^{\infty}j\frac{w^q(Q)^{\frac{1}{q}-\frac{\kappa}{p}}}{w^q(2^{j+1}Q)^{\frac{1}{q}-\frac{\kappa}{p}}}\\
&\leq C\sum\limits_{|\gamma|=m-1}\|D^{\gamma}A\|_*\|f\|_{L^{p,\kappa}(w^p,w^q)}\sum\limits_{j=1}^{\infty}\frac{j}{{(D_1^{j+1})}^{\frac{1}{q}-\frac{\kappa}{p}}}
\leq C\sum\limits_{|\gamma|=m-1}\|D^{\gamma}A\|_*\|f\|_{L^{p,\kappa}(w^p,w^q)}
\end{align*}
where $D_1>1$ is the reverse doubling constant.
Taking supremum over all cubes in $\mathbb{R}^n$ on both sides of the above inequality, we complete the proof of inequality (\ref{eq3.5}) of Theorem \ref{th3.3}.\\
It is not difficult to see that inequality (\ref{eq3.6}) is easy to get from (\ref{eq3.5}) and (\ref{eq5.1}).\\
{\bf{The proof of Theorem \ref{th3.4}}}\qquad We consider (\ref{eq3.7}) firstly. Let $Q$ be the same as in the proof of (\ref{eq3.5}) and $\bar{Q}=2Q$, we decompose $f$ according to $\bar{Q}$: $f=f\chi_{\bar{Q}}+f\chi_{(\bar{Q})^c}:=f_1+f_2$. Thus we have
\begin{align*}
\|[A,T_{\Omega}]f\|_{L^{p,\kappa}(w)}&\leq \frac{1}{w(Q)^{\kappa/p}}\Big(\int_{Q}|[A,T_{\Omega}]f_1(y)|^pw(y)dy\Big)^{\frac{1}{p}}\\
&\qquad+\frac{1}{w(Q)^{\kappa/p}}\Big(\int_{Q}|[A,T_{\Omega}]f_2(y)|^pw(y)dy\Big)^{\frac{1}{p}}=I+II
\end{align*}
From Theorem \ref{lem10} (a), the $L^p(w)$ boundedness of $[A,T_{\Omega}]$ and Lemma \ref{lem12} (1) that $w\in \Delta_2$, we get
\begin{align*}
&I\leq \frac{1}{w(Q)^{\kappa/p}}\|[A,T_{\Omega}]f_1\|_{L^p(w)}\leq \frac{C}{w(Q)^{\kappa/p}}\|A\|_*\|f_1\|_{L^p(w)}
\end{align*}
\begin{align*}
&=C\|A\|_*\|f\|_{L^{p,\kappa}(w)}\frac{w(\bar{Q})^{\kappa/p}}{w(Q)^{\kappa/p}}\leq C\|A\|_*\|f\|_{L^{p,\kappa}(w)}
\end{align*}
For $|[A,T_{\Omega}]f_2(y)|$, using H$\ddot{o}$lder inequality, we obtain
\begin{align*}
&|[A,T_{\Omega}]f_2(y)|\leq \sum\limits_{j=1}^{\infty}\int_{2^{j+1}Q\setminus2^jQ}\frac{|\Omega(y-z)|}{|y-z|^n}|A(y)-A(z)||f(z)|dz\\
&\leq C\sum\limits_{j=1}^{\infty}\frac{1}{|2^jQ|}(\int_{2^{j+1}Q}|\Omega(y-z)|^sdz)^{\frac{1}{s}}(\int_{2^{j+1}Q}|A(y)-A(z)|^{s'}|f(z)|^{s'}dz)^{\frac{1}{s'}}\\
&\leq C\|\Omega\|_{L^s}\sum\limits_{j=1}^{\infty}\frac{|2^{j+2}Q|^{\frac{1}{s}}}{|2^jQ|}\Big|A(y)-m_{2^{j+1}Q,w^{-\frac{s'}{p-s'}}}(A)\Big|\Big(\int_{2^{j+1}Q}|f(z)|^{s'}dz\Big)^{\frac{1}{s'}}\\
&+C\|\Omega\|_{L^s}\sum\limits_{j=1}^{\infty}\frac{|2^{j+2}Q|^{\frac{1}{s}}}{|2^jQ|}\Big(\int_{2^{j+1}Q}|m_{2^{j+1}Q,w^{-\frac{s'}{p-s'}}}(A)-A(z)|^{s'}|f(z)|^{s'}dz\Big)^{\frac{1}{s'}}
&=II_1(y)+II_2
\end{align*}
We estimate $II_1(y)$ and $II_2$ respectively. By H$\ddot{o}$lder inequality and the condition that $w\in A_{p/s'}$, we have
\begin{align*}
&\frac{1}{w(Q)^{\kappa/p}}\Big(\int_{Q}II_1(y)^pw(y)dy\Big)^{\frac{1}{p}}\\
&=C\frac{\|\Omega\|_{L^s}}{w(Q)^{\kappa/p}}\sum\limits_{j=1}^{\infty}\frac{|2^{j+2}Q|^{\frac{1}{s}}}{|2^jQ|}\Big(\int_Q|A(y)-m_{2^{j+1}Q,w^{-\frac{s'}{p-s'}}}(A)|^p(\int_{2^{j+1}Q}|f(z)|^{s'}dz)^{\frac{p}{s'}}w(y)dy\Big)^{\frac{1}{p}}\\
\end{align*}
\begin{align*}
&\leq C\frac{\|f\|_{L^{p,\kappa}(w)}}{w(Q)^{\kappa/p}}\sum\limits_{j=1}^{\infty}\frac{|2^{j+2}Q|^{\frac{1}{s}}}{|2^jQ|}w(2^{j+1}Q)^{\frac{\kappa}{p}}w^{-\frac{s'}{p-s'}}(2^{j+1}Q)^{\frac{p-s'}{ps'}}\cdot\\
&\qquad\Big(\int_Q|A(y)-m_{2^{j+1}Q,w^{-\frac{s'}{p-s'}}}(A)|^pw(y)dy\Big)^{\frac{1}{p}}\\
&\leq C\frac{\|f\|_{L^{p,\kappa}(w)}}{w(Q)^{\kappa/p}}\sum\limits_{j=1}^{\infty}\frac{1}{w^(2^{j+1}Q)^{\frac{1-\kappa}{p}}}\Big(\int_Q|A(y)-m_{2^{j+1}Q,w^{-\frac{s'}{p-s'}}}(A)|^pw(y)dy\Big)^{\frac{1}{p}}\\
\end{align*}
We calculate the part of including $m_{2^{j+1}Q,w^{-\frac{s'}{p-s'}}}(A)$ as follows:
\begin{align*}
&\Big(\int_Q|A(y)-m_{2^{j+1}Q,w^{-\frac{s'}{p-s'}}}(A)|^pw(y)dy\Big)^{\frac{1}{p}}\\
&\leq \Big(\int_Q|A(y)-m_{Q,w}(A)|^pw(y)dy\Big)^{\frac{1}{p}}+\Big|m_{Q,w}(A)-m_{2^{j+1}Q,w^{-\frac{s'}{p-s'}}}(A)\Big|w(Q)^{\frac{1}{p}}\\
&=III+IV
\end{align*}
For term $III$, notice that $w\in A_{p/s'}\subset A_{\infty}$, thus from Lemma \ref{lem13}, we get
\begin{align*}
III\leq w(Q)^{\frac{1}{p}}\|A\|_*
\end{align*}
Then we estimate $IV$. By Lemma \ref{lem9} and Lemma \ref{lem13}, we have
\begin{align*}
&\Big|m_{Q,w}(A)-m_{2^{j+1}Q,w^{-\frac{s'}{p-s'}}}(A)\Big|\leq|m_{Q,w}(A)-m_Q(A)|+|m_Q(A)-m_{2^{j+1}Q}(A)|\\
&\qquad+\Big|m_{2^{j+1}Q}(A)-m_{2^{j+1}Q,w^{-\frac{s'}{p-s'}}}(A)\Big|\\
&\leq \frac{1}{w(Q)}\int_Q|A(t)-m_Q(A)|w(t)dt+2^n(j+1)\|A\|_*\\
&\qquad+\frac{1}{w^{-\frac{s'}{p-s'}}(2^{j+1}Q)}\int_{2^{j+1}Q}|A(t)-m_{2^{j+1}Q}(A)|w^{-\frac{s'}{p-s'}}(t)dt\\
&\leq C(j+1)\|A\|_*
\end{align*}
So,
\begin{align*}
IV\leq C(j+1)\|A\|_*w(Q)^{\frac{1}{p}}
\end{align*}
As a result,
\begin{align*}
&\frac{1}{w(Q)^{\kappa/p}}\Big(\int_{Q}II_1(y)^pw(y)dy\Big)^{\frac{1}{p}}\\
&\leq C\|A\|_*\|f\|_{L^{p,\kappa}(w)}\sum\limits_{j+1}^{\infty}(j+1)\frac{w(Q)^{\frac{1-\kappa}{p}}}{w(2^{j+1}Q)^{\frac{1-\kappa}{p}}}\\
&\leq C\|A\|_*\|f\|_{L^{p,\kappa}(w)}
\end{align*}
For $II_2$, by H$\ddot{o}$lder inequality and $w\in A_{p/s'}$, we get
\begin{align*}
&II_2\leq C\|\Omega\|_{L^s}\sum\limits_{j=1}^{\infty}\frac{|2^{j+2}Q|^{\frac{1}{s}}}{|2^jQ|}(\int_{2^{j+1}Q}|f(z)|^pw(z)dz)^{\frac{1}{p}}\\
&\qquad\Big(\int_{2^{j+1}Q}\Big|m_{2^{j+1}Q,w^{-\frac{s'}{p-s'}}}(A)-A(z)\Big|^{\frac{ps'}{p-s'}}w(z)^{-\frac{s'}{p-s'}}dz\Big)^{\frac{p-s'}{ps'}}\\
&\leq C\|\Omega\|_{L^s}\|A\|_*\|f\|_{L^{p,\kappa}(w)}\sum\limits_{j=1}^{\infty}\frac{|2^{j+2}Q|^{\frac{1}{s}}}{|2^jQ|}w(2^{j+1}Q)^{\frac{\kappa}{p}}w^{-\frac{s'}{p-s'}}(2^{j+1}Q)^{\frac{p-s'}{ps'}}\\
&\leq C\|A\|_*\|f\|_{L^{p,\kappa}(w)}\sum\limits_{j=1}^{\infty}\frac{1}{w(2^{j+1}Q)^{\frac{1-\kappa}{p}}}
\end{align*}
Therefore,
\begin{align*}
&\frac{1}{w(Q)^{\kappa/p}}\Big(\int_{Q}II_2^pw(y)dy\Big)^{\frac{1}{p}}\\
&\leq C\|A\|_*\|f\|_{L^{p,\kappa}(w)}\sum\limits_{j=1}^{\infty}\frac{w(Q)^{\frac{1-\kappa}{p}}}{w(2^{j+1}Q)^{\frac{1-\kappa}{p}}}\\
&\leq C\|A\|_*\|f\|_{L^{p,\kappa}(w)}
\end{align*}
So far, we have completed the proof of (\ref{eq3.7}).\\
The inequality (\ref{eq3.8}) can be immediately obtained from (\ref{eq5.1}) and (\ref{eq3.7}).

{\bf{The proof of Theorem \ref{th3.5}}}\qquad As before, we prove (\ref{eq3.9}) at first. Assume $Q$ be the same as in the proof of (\ref{eq3.5}) and $\bar{Q}=2Q$, set
\begin{align*}
A_{\bar{Q}}(y)=A(y)-m_{\bar{Q}}(\nabla A)y
\end{align*}
We also decompose $f$ according to $\bar{Q}$ : $f=f\chi_{\bar{Q}}+f\chi_{(\bar{Q})^c}:=f_1+f_2$. Then we get
\begin{align*}
\|\tilde{T}_{\Omega}^Af\|_{L^{p,\kappa}(w)}&\leq \frac{1}{w(Q)^{\kappa/p}}\Big(\int_{Q}|\tilde{T}_{\Omega}^Af_1(y)|^pw(y)dy\Big)^{\frac{1}{p}}\\
&\qquad+\frac{1}{w(Q)^{\kappa/p}}\Big(\int_{Q}|\tilde{T}_{\Omega}^Af_2(y)|^pw(y)dy\Big)^{\frac{1}{p}}=I+II
\end{align*}
For the first term $I$, Theorem \ref{lem11} and Lemma \ref{lem12} (1) implies
\begin{align*}
&I\leq \frac{1}{w(Q)^{\kappa/p}}\|\tilde{T}_{\Omega}^Af_1\|_{L^p(w)}\leq \frac{C}{w(Q)^{\kappa/p}}\|\Omega\|_{\infty}\|\nabla A\|_*\|f_1\|_{L^p(w)}\\
&\leq C\|\nabla A\|_*\|f\|_{L^{p,\kappa}(w)}\frac{w(\bar{Q})^{\kappa/p}}{w(Q)^{\kappa/p}}\leq C\|\nabla A\|_*\|f\|_{L^{p,\kappa}(w)}
\end{align*}
We will omit the proof for the term $II$ as it is similar to and easier than the part of $II$ in the proof of (\ref{eq3.5}), except under the condition that $w\in A_p$, $\Omega\in L^{\infty}(S^{n-1})$, $m=2$ and $f\in L^{p,\kappa}(w)$. For inequality (\ref{eq3.10}), it can be easily proved by (\ref{eq3.9}) and (\ref{eq5.1}). Thus, we complete the proof of Theorem \ref{th3.5}.

\vspace{1cm} \noindent {\bf {References}} \small

\begin{itemize}
\item [{[1]}] B\,Bajsanski, R\,Coifman. \emph{On singular integrals}, Proc Symp Pure Math Amer Math Soc Providence, R I, 1967, 38(10): 1-17.
\item [{[2]}] J\,Cohen. \emph{A sharp estimate for a multilinear singular integral in $R^n$}, Indiana Univ Math J, 1981, 30(5): 693-702.
\item [{[3]}] J\,Cohen, J\,Gosselin. \emph{A BMO estimate for a multilinear singular integrals}, Ill J Math, 1986, 30(3): 445-464.
\item [{[4]}] S\,Chanillo. \emph{A note on commutators}, Indiana Univ Math J, 1982, 31: 7-16.
\item [{[5]}] Y\,Ding. \emph{A note on multilinear fractional integrals with rough kernel}, Advan Math, 2001, 30(3): 238-246.
\item [{[6]}] Y\,Ding. \emph{A two-weight estimate for a class of fractional integral operators with rough kernel}, Int J Math Math Sci, 2005, 2005(12): 1835-1842.
\item [{[7]}] Y\,Ding, S\,Z\,Lu. \emph{Weighted norm inequalities for fractional integral operators with rough kernel}, Can J Math, 1998, 50(1): 29-39.
\item [{[8]}] Y\,Ding, S\,Z\,Lu. \emph{Weighted boundedness for a class of rough multilinear operators}, Acta Math Sin-English Series, 2001, 17(3): 517-526.
\item [{[9]}] A\,DeVore Ronald, C\,Sharpley Robert. \emph{Maximal functions measuring smoothness}, Mem Amer Math Soc, 1984, 47(293): 115.
\item [{[10]}] L\,Grafakos. \emph{Classical and Modern Fourier Analysis} (Pearson Education), Inc Upper Saddle River, New Jersey, 2004.
\item [{[11]}] H\,Y\,Han, S\,Z\,Lu. \emph{Lipschitz estimates for commutators of fractional integral operators with rough kernel}, J Beijing Normal Univ (Natural Science), 2006, 42(1): 11-18.
\item [{[12]}] S\,Hofmann. \emph{On certain nonstandard Calder\'{o}n-Zygmund operators}, Stud Math, 1994, 109(2): 105-131.
\item [{[13]}] Y\,L\,Jiao. \emph{A note on the multilinear singular integral operators}, Ann Theo Appl, 2004, 20(4): 373-382.
\item [{[14]}] Y\,Komori, S\,Shirai. \emph{Weighted Morrey spaces and a singular integral operator}, Math Nachr, 2009, 282(2): 219-231.
\item [{[15]}] S\,Z\,Lu, Y\,Ding, D\,Y\,Yan. \emph{Singular integrals and related topics}, World Scientific, 2007.
\item [{[16]}] S\,Z\,Lu, H\,X\,Wu, P\,Zhang. \emph{Multilinear singular integrals with rough kernel}, Acta Math Sin-English Series, 2003, 19(1): 51-62.
\item [{[17]}] S\,Z\,Lu, P\,Zhang. \emph{Lipschitz estimates for generalized commutators of fractional integrals with rough kernel}, Math Nachr, 2003, 252(1): 70-85.
\item [{[18]}] B\,Muckenhoupt, R\,Wheeden. \emph{Weighted norm inequalities for fractional integrals}, Trans Amer Math Soc, 1974, 192: 261-274.
\item [{[19]}] B\,Muckenhoupt, R\,Wheeden. \emph{Weighted bounded mean oscillation and the Hilbert transform}, Stud Math, 1976, 54: 221-237.
\item [{[20]}] M\,Paluszynski. \emph{Characterization of the Besov space via the commutator operator of Coifman, Rochberg and Weiss}, Indiana Univ Math J, 1995, 44: 1-17.
\item [{[21]}] Y\,L\,Shi, X\,X\,Tao. \emph{Some multi-sublinear operators on generalized Morrey spaces with non-doubling measures}, J Korean Math Soc, 2012, 49(5): 907-925.
\item [{[22]}] E\,M\,Stein. \emph{Harmonic Analysis: Real-variable Methods, Orthogonality, and Oscillatory Integrals}, Princeton Univ Press, Princeton, NJ, 1993.
\item [{[23]}] X\,X\,Tao, Y\,L\,Shi. \emph{Multilinear Riesz potential operators on Herz-type spaces and generalized Morrey spaces}, Hokkaido Math J, 2009, 38(3): 635-662.
\item [{[24]}] X\,X\,Tao, Y\,L\,Shi. \emph{Multilinear commutators of calder\'{o}n-zygmund operator on
    $\lambda$-central Morrey spaces}, Advan Math (China), 2011, 40(1): 47-59.
\item [{[25]}] X\,X\,Tao, Y\,L\,Shi, S\,Y\,Zhang. \emph{Boundedness of multilinear Riesz potential operators on product of Morrey spaces and Herz-Morrey spaces}, Acta Math Sin-Chinese Series, 2009, 52(3): 535-548.
\item [{[26]}] X\,X\,Tao, Y\,L\,Shi, T\,T\,Zheng. \emph{Multilinear Riesz Potential on Morrey-Herz spaces with non-doubling measures}, J Inequal Appl, 2010, 2010(Article ID 731016): 21 pages.
\item [{[27]}] X\,X\,Tao, H\,H\,Zhang. \emph{On the boundedness of multilinear operators on weighted Herz-Morrey spaces}, Taiwan J Math, 2011, 15(4): 1527-1543.
\item [{[28]}] X\,X\,Tao, T\,T\,Zheng. \emph{Multilinear commutators of fractional integrals over Morrey spaces with non-doubling measures}, Nodea-Nonlinear Differ Equ Ap, 2011, 18(3): 287-308.
\item [{[22]}] A\,Torchinsky. \emph{Real variable methods in Harmonic Analysis}, Academic Press, San Diego, 1986.
\item [{[23]}] H\,Wang. \emph{The boundedness of some operators with rough kernel on the weighted Morrey spaces}, arXiv:1011.5763v1[math.CA].
\item [{[24]}] Q\,Wu, D\,C\,Yang. \emph{On fractional multilinear singular integrals}, Math Nachr, 2002, 239-240(1): 215-235.
\item [{[25]}] Y\,P\,Wu, X\,X\,Tao. \emph{Estimate for a class of multilinear fractional integral operator with rough kernel}, Ann Theo Appl, 2010, 26(4): 359-366.

\end{itemize}\vskip 10mm Department of
Mathematics, Hangzhou Normal University, Hangzhou 310023,
China\\
\indent Email: amyhesha@163.com\\
Department of Mathematics, Zhejiang University of Science and Technology, Hangzhou 310023,
China\\
\indent Email: xxtao@zust.edu.cn

\end{document}